\documentclass[12pt,reqno,a4paper]{amsart} 
\usepackage{amsmath,amssymb}
\usepackage{fullpage}
\usepackage{enumerate}
\usepackage{url}

\usepackage[english]{babel}

\newcommand{\dx}{\mathrm{d}} 
\newcommand{\eps}{\varepsilon}
\newcommand{\sumast}{\sideset{}{^*}\sum}
\newcommand{\sumprime}{\sideset{}{^\prime}\sum}
\newcommand{\singseries}{\mathfrak{S}}

\newcommand{\Odip}[2]{\mathcal{O}_{#1}\!\left(#2\right)\mathchoice{\!}{}{}{}}
\newcommand{\Oprimedip}[2]{\mathcal{O}^{\prime}_{#1}\!\left(#2\right)\mathchoice{\!}{}{}{}}

\newcommand{\Odi}[1]{\Odip{}{#1}}
\newcommand{\Oprimedi}[1]{\Oprimedip{}{#1}}

\newenvironment{Proof}[1][Proof]{\par\noindent\textbf{#1.}~}
  {\hfill$\square$\smallskip\par}

\newtheorem{Lemma}{Lemma} 

 \newtheoremstyle{Nonumtheorems}
  {10pt}
  {6pt}
  {\itshape}
  {}
  {\bfseries}
  {.}
  {.5em}
  {\thmname{#1}\thmnote{ (#3)}}

\theoremstyle{Nonumtheorems}
\newtheorem{Nonumthm}{Theorem}

\begin{document}
\allowdisplaybreaks
\title{Sums of many primes}
\author[A.~Languasco \& A.~Zaccagnini]{Alessandro Languasco \& Alessandro Zaccagnini}
\begin{abstract}
Assuming that the Generalized Riemann Hypothesis (GRH) holds,
we prove an explicit formula for the number of representations of 
an integer as a sum of $k\geq 5$ primes.
Our error terms in such a formula improve by some logarithmic factors an analogous result by  Friedlander-Goldston \cite{FriedlanderG1997}.
\end{abstract}
\maketitle
\section{Introduction}

Let $k \ge 2$ be a fixed integer and set
\begin{equation}
\label{Rk-def}
  R_k(n)
  =
  \sum_{n_1 + \cdots + n_k = n}
    \Lambda(n_1) \cdots \Lambda(n_k)
\end{equation}
and
\begin{equation}
\label{singseries-def}
  \singseries_k(n)
  =
  \sum_{q = 1}^{\infty}
    \frac{\mu(q)^k}{\phi(q)^k} c_q(-n)
  =
  \prod_{p \mid n}
    \Bigl( 1 - \Bigl( \frac{-1}{p - 1} \Bigr)^{k-1}
    \Bigr)
  \prod_{p \nmid n}
    \Bigl( 1 - \Bigl( \frac{-1}{p - 1} \Bigr)^k
    \Bigr)
\end{equation}
where $c_q$ is the Ramanujan sum defined as
\begin{equation}
\label{cq-def}
    c_q(m) = \sumast_{a = 1}^q e \Bigl( m \frac aq \Bigr).
\end{equation}
Moreover let $\chi \bmod{q}$ be a Dirichlet character and
\begin{equation}
\label{c-chi-def}
  c_{\chi}(m)
  =
  \sum_{a = 1}^q \chi(a) e \Bigl( m \frac aq \Bigr). 
\end{equation}
Since we deal with the case $k = 2$ in \cite{LanguascoZaccagnini2010},
here we assume that $k \ge 3$ throughout for simplicity of statement.
We just notice that in \cite{LanguascoZaccagnini2010} there is an
average of $R_k(n)$ over $n$ and the natural hypothesis to make is RH,
whereas here and in Friedlander-Goldston \cite{FriedlanderG1997} there
is no such average and the natural hypothesis is GRH.
In both cases we are interested into a formula which is ``explicit''
in the sense that it has the expected main term, a secondary main term
depending on the zeros of the $L$ functions (or just the zeta function
when $k = 2$), and an error term of smaller order of magnitude.

We have, for $k\geq 5$, the following explicit formula for  $R_k(n)$.
\begin{Nonumthm}
\label{Explicit-formula-Theorem}
Let $k\geq 5$ be a fixed integer.
Assume that the Generalized Riemann Hypothesis (GRH) holds 
for every Dirichlet series $L(s,\chi)$, for every
$\chi \bmod{q}$.
Then, for every sufficiently large integer $n$,  we have that
\begin{align}
\notag
  R_k(n)
  &=
   \frac{n^{k - 1}}{(k - 1)!}
  \singseries_k(n)
- k
  \sum_{q = 1}^{\infty}
    \frac{\mu(q)^{k - 1}}{\phi(q)^k}
    \sum_{\chi \bmod q}
      c_{\chi}(-n) \tau(\overline{\chi})
      \sum_{\rho}
        \frac{n^{\rho+k-2}}{\rho(\rho + 1) \cdots (\rho + k - 2)} 
        \\
\label{expl-form-Rk}
  &\qquad+
  \Odi{n^{k - 7/4}\log ^{k-1} n}
\end{align}
where $R_k(n),   \singseries_k(n), c_{\chi}(n)$ are respectively
defined in \eqref{Rk-def}-\eqref{singseries-def} and \eqref{c-chi-def},
$\tau(\chi)$ is the Gauss sum
and  $\rho=1/2+i\gamma$ runs over
the non-trivial zeros of $L(s,\chi)$,
for every $\chi \bmod{q}$.
%
For $k\ge 6$, in the last error term we can replace $7/4$ by $2$.
For $k\ge 7$, in the last error term we can also replace $\log ^{k-1} n$ by $\log ^{2} n$.
\end{Nonumthm} 

The condition $k\geq 5$ essentially arises in two points. The first
one is the evaluation of the secondary main term in \eqref{expl-form-Rk} 
(see the error term in \eqref{M1-eval} and the remark at the end
of \S \ref{secondary-main-term}) while the second one is in the error term 
estimates (see \S \ref{error-terms}).

This result should be compared with Proposition~1
of Friedlander-Goldston \cite{FriedlanderG1997}. 
They have a more involved but
equivalent form of the secondary main term and worse
estimates for the error term.
In principle, both here and in \cite{FriedlanderG1997} one could give
a statement with the sum over $q$ in the ``secondary main term'' in
the right hand side of \eqref{expl-form-Rk} restricted to
$q \le n^{1/2} / 2$, and assume only that the GRH holds for the $L$
functions associated to characters modulo these values of $q$.
For the details, see the remark at the end of
\S\ref{secondary-main-term}.

The improvement given here is due to the fact that we use the version
of the circle method introduced by Hardy and Littlewood in
\cite{HardyL1923} and used also by Linnik in
\cite{Linnik1946,Linnik1952}, involving series rather than truncated
sums: it is essentially equivalent to, but slightly sharper than, the
usual approach with truncated sums.

\textbf{Acknowledgments.} 
We would like to thank the organizers of the conference 
``Number Theory and its applications, 
An International Conference Dedicated to 
K\'alm\'an Gy\H{o}ry, Attila Peth\H{o}, J\'anos Pintz, Andr\'as S\'ark\"ozy'',
for the excellent working environment provided.

\section{Lemmas}

We will  use the original Hardy and Littlewood \cite{HardyL1923} circle 
method setting, \emph{i.e.},
the weighted exponential sum
\begin{equation}
\label{tildeS-def}
\widetilde{S}(\alpha)
=
\sum_{n=1}^{\infty} 
\Lambda(n) e^{-n/N}
e(n\alpha),
\end{equation}
where $e(x)=\exp(2\pi i x)$,
since it lets us avoid the use of Gallagher's Lemma 
(Lemma 1 of \cite{Gallagher1970}) and hence it gives slightly sharper
results in this conditional case: see Lemma \ref{LP-Lemma} below.

Let $1\leq Q\leq N$ be a parameter to be chosen later.
We will consider the set of the Farey fractions of level $Q$ 
\[
\left\{
\frac{a}{q}: 1\leq q\leq Q, 0\leq a\leq q, (a,q)=1
\right\}.
\]
Let $a'/q'<  a/q < a''/q''$ be three consecutive
Farey fractions, 
\[
{\mathcal M}_{q,a}=
\left(
\frac{a+a'}{q+q'},\frac{a+a''}{q+q''}
\right]
\quad
\textrm{if}\quad
\frac{a}{q}\neq \frac11
\] 
and
${\mathcal M}_{1,1}=(1-1/(Q+1),1+1/(Q+1)]$
be the Farey arcs centered at $a/q$. These intervals are disjoint
and their union is $(1/(Q+1),1+1/(Q+1)]$.
Moreover, let
\begin{equation}
\label{Farey-arc}
\xi_{q,a}=
\left(
\frac{-1}{q(q+q')},\frac{1}{q(q+q'')}
\right]
\end{equation}
and
$\xi_{1,1}=(-1/(Q+1),1/(Q+1)]$
be the Farey arcs re-centered at the origin.
In the following we also use the relation
\[
\Bigl(\frac{-1}{2qQ},\frac{1}{2qQ}\Bigr)
\subseteq
\xi_{q,a}
\subseteq
\Bigl(\frac{-1}{qQ},\frac{1}{qQ}\Bigr).
\]
Let
\begin{equation}
\label{z-def}
  z
  =
  N^{-1} - 2 \pi i \eta
\end{equation}
for $\eta\in \xi_{q,a}$, and
\begin{equation}
\notag
  V(\eta)
  =
  \sum_{m = 1}^{\infty} e^{-m / N} e(m \eta)
  =
  \sum_{m = 1}^{\infty} e^{-m z}
  =
  \frac1{e^z - 1}.
\end{equation}

\begin{Lemma}
\label{V-behaviour}
If $z$ satisfies \eqref{z-def} then $V(\eta) = z^{-1} + \Odi{1}$.
\end{Lemma}

\begin{Proof}
We recall that the function $w / (e^w - 1)$ has a power-series
expansion with radius of convergence $2 \pi$ (see for example Apostol
\cite{Apostol1976}, page 264).
In particular, uniformly for $|w| \le 4 < 2 \pi$ we have
$w / (e^w - 1) = 1 + \Odi{|w|}$.
Since $z$ satisfies \eqref{z-def} we have $|z| \le 4$ and the result
follows.
\end{Proof}

Combining Lemma \ref{V-behaviour} and the inequality
\begin{equation}
\label{z-estim}
\vert z\vert ^{-1} \ll \min \bigl(N, \vert \eta \vert^{-1}\bigr),
\end{equation}
we also have
\begin{equation}
\label{V-estim}
\vert V(\eta) \vert \ll \vert z\vert ^{-1}+1 \ll \min \bigl(N, \vert \eta \vert^{-1}\bigr).
\end{equation}

We will use the approximation
\begin{equation}
\label{approx-with-V-1}
  \widetilde{S} \Bigl( \frac aq + \eta \Bigr)
  =
  \frac{\mu(q)}{\phi(q)}V(\eta) 
  +
  \widetilde{R}(\eta; q, a, V),
\end{equation}
and
\begin{align}
\label{approx-with-V-2}
\widetilde{R}(\eta; q, a, V)
=
  \frac1{\phi(q)}
  \sum_{\chi \bmod q}
    \chi(a) \tau(\overline{\chi})
	W(\chi,\eta,V)
+
  \Odi{(\log(q N))^2},
\end{align}
where 
\[
W(\chi,\eta, V)
=   
      \sum_{\ell = 1}^{\infty}
        \Lambda(\ell) \chi(\ell) e^{-\ell / N} e(\ell \eta)
      -
        \delta(\chi)V(\eta),
\]
$\delta(\chi)=1$ if $\chi=\chi_{0} \bmod{q}$ and $0$ otherwise.
Recalling \eqref{z-def}, by Lemma \ref{V-behaviour} we can also write
\[
  \widetilde{S} \Bigl( \frac aq + \eta \Bigr)
  =
  \frac{\mu(q)}{\phi(q)z}
  +
  \widetilde{R}(\eta; q, a, z) +\Odi{\frac{1}{\phi(q)}}
\]
and
\begin{equation}
\notag
  \widetilde{R}(\eta; q, a, z)
 =
\frac1{\phi(q)}
  \sum_{\chi \bmod q}
    \chi(a) \tau(\overline{\chi})
W(\chi,\eta,z)  
+
  \Odi{(\log(q N))^2},
\end{equation}
where
\[
W(\chi,\eta,z)
= 
	\sum_{\ell = 1}^{\infty}
        \Lambda(\ell) \chi(\ell) e^{-\ell / N} e(\ell \eta)
      -
        \frac{\delta(\chi)}{z}.
\]
Summing up we have 
\begin{equation}
\label{approx-with-z-2}
  \widetilde{S} \Bigl( \frac aq + \eta \Bigr)
  =
  \frac{\mu(q)}{\phi(q)z}
  +
  \frac1{\phi(q)}
  \sum_{\chi \bmod q}
    \chi(a) \tau(\overline{\chi})
  W(\chi,\eta,z)
  +
  \Odi{(\log(q N))^2}.
\end{equation}

Recalling \eqref{tildeS-def}, the
first ingredient we need is the following explicit formula
which slightly sharpens what Linnik \cite{Linnik1946}
(see also  eq.~(4.1) of \cite{Linnik1952}) proved.
\begin{Lemma} 
\label{Linnik-lemma}
If $\chi$ is a character ${}\bmod q$ and GRH holds for $L(s, \chi)$
then 
\begin{equation}
\label{expl-form}
W(\chi,\eta,z)
=
- \sum_{\rho}z^{-\rho}\Gamma(\rho) 
+
E(q, N)
\end{equation}
where $\rho=\beta+i\gamma$ runs over the non-trivial zeros of
$L(s,\chi)$ and
\begin{equation}
\label{expl-form-err-term-bis}
E(q, N)
  \ll
  \begin{cases}
    1               + \log^2 q & \text{if $\chi$ is a primitive character,} \\
    (\log N) (\log q) + \log^2 q & \text{if $\chi$ is not primitive.}
  \end{cases}
\end{equation}
\end{Lemma}
\begin{Proof}
We recall that $\delta(\chi) = 1$ if $\chi = \chi_{0} \bmod{q}$ and
$0$ otherwise.
Let
\[
  \Sigma(N, \chi, \eta)
  =
  \sum_{\ell = 1}^{\infty} \Lambda(\ell) \chi(\ell) e^{-\ell / N} e(\ell \eta)
  =
  W(\chi, \eta, z)
  +
  \frac{\delta(\chi)}z.
\]
We notice that if $\chi \bmod q$ is induced by $\chi_1 \bmod q_1$ then
\begin{align*}
  \bigl| \Sigma(N, \chi, \eta) - \Sigma(N, \chi_1, \eta) \bigr|
  &\le
  \sum_{\substack{\ell \ge 1 \\ (\ell, q) > 1 \\ (\ell, q_1) = 1}} \Lambda(\ell) e^{-\ell / N}
  \ll
  \log q \log N.
\end{align*}
We now assume that $\chi \bmod q$ is a primitive character and let
$\alpha = 3 / 4$.
Following the proof of Lemma 4 in Hardy and Littlewood \cite{HardyL1923}
and \S4 in Linnik \cite{Linnik1946}, we have that
\begin{equation}
\label{Mellin-bis}
  W(\chi, \eta, z)
  =
  - 
  \sum_{\rho}z^{-\rho} \Gamma(\rho) 
  +
  C(\chi)
  -
  \frac{1}{2\pi i}
  \int_{(-\alpha)} 
  \frac{L'}{L}(w, \chi) \Gamma(w) z^{-w} \, \dx w,
\end{equation}
where $C(\chi)$ is a term that depends only on the character $\chi$.
In order to estimate the integral in \eqref{Mellin-bis} we need the
inequality
\begin{equation}
\label{L'/L-bound}
  \left\vert \frac{L'}L \Bigl( -\frac34 + i t, \chi \Bigr)
  \right\vert
  \ll
  \log \bigl( q (\vert t \vert + 2 ) \bigr).
\end{equation}
This follows from equations (1) and (4) of \S16 of Davenport
\cite{Davenport2000} since the latter reads
\[
  \frac{L'}{L}(w, \chi)
  =
  \sumprime_\rho \frac1{w - \rho}
  +
  \Odi{\log( q (\vert t \vert + 2) )},
\]
where the dash means that the sum is restricted to those zeros
$\rho = \beta + i \gamma$ with $|t - \gamma| < 1$, while the former
implies that the number of such summands is
$\ll \log( q (\vert t \vert + 2) )$.
Finally, it is obvious that each summand is $\ll 1$ on the line of
integration $w = -\alpha + i t$.

We notice that $|z^{-w}|= \vert z \vert^{\alpha} \exp(t \arg(z))$
where $\vert \arg(z) \vert \le \frac12 \pi$.
Furthermore the Stirling formula implies that
$\Gamma(w) \ll \vert t \vert^{-\alpha-1/2} \exp(-\frac{\pi}{2}\vert t \vert)$.
Hence
\begin{align*}
  \int_{(-\alpha)}
    \frac{L'}{L}(w,\chi) \Gamma(w) z^{-w} \, \dx w
  &\ll
  |z|^{\alpha}
  \int_0^1 \log(q (t + 2)) \, \dx t \\
  &\qquad
  +
  |z|^{\alpha}
  \int_1^{\infty}
    \log(q (t + 2)) t^{-\alpha-1/2}
    \exp\Bigl( (\arg(z) - \frac{\pi}{2}) t \Bigr) \, \dx t \\
  &\ll
  |z|^{\alpha} (1 + \log q)
  +
  |z|^{\alpha}
  \int_1^{\infty}
    \log(q (t + 2)) t^{-\alpha-1/2} \, \dx t \\
  &\ll
  |z|^{\alpha} (1 + \log q).
\end{align*}
This is $\ll 1 + \log(q)$ as stated since $z \ll 1$ by \eqref{z-def}
and $\alpha$ is fixed.
Finally, we have to deal with the term $C(\chi)$ in
\eqref{Mellin-bis}.
We recall the notation of \S19 of Davenport \cite{Davenport2000}: if
$\chi$ is odd then $b(\chi)$ denotes $(L' / L)(0, \chi)$, whereas if
$\chi$ is even it is the constant term in the Laurent expansion of
$(L' / L)(w, \chi)$ around zero.
In other words, $(L' / L)(w, \chi) = w^{-1} + b(\chi) + \Odi{w}$.

If $\chi$ is odd $C(\chi)$ is simply $-(L' / L)(0, \chi) = -b(\chi)$
since $L(0, \chi) \ne 0$.
If $\chi$ is even then $L(w, \chi)$ has a simple zero at $0$ and
therefore $-(L' / L)(w, \chi) \Gamma(w) z^{-w}$ has a double pole at
$w = 0$ with residue $C(\chi) = \log(z) - b(\chi) - \Gamma'(1)$.
Arguing as on pages 118--119 of Davenport \cite{Davenport2000}, we see
that
\[
  b(\chi)
  =
  -\sum _{\rho} \left(\frac{1}{\rho} + \frac{1}{2-\rho}\right)
  +
  \Odi{1}
  =
  -
  \sum_{\vert \gamma \vert < 1} \frac{1}{\rho}
  +
  \Odi{\log q}
  \ll
  \log^2 q.
\]
Finally, $\log(z) \ll 1$ since $z$ satisfies \eqref{z-def}.
\end{Proof}
\begin{Lemma} 
\label{incond-mean-square}
Let $N$ be a sufficiently large integer, $Q\leq N$ and $z$ be as in \eqref{z-def}
with $\eta \in \xi_{q,a}$. 
We have 
\[
 \sum_{q = 1}^Q 
    \sumast_{a = 1}^q  
    \int_{\xi_{{q,a}}}
\Bigl\vert
\frac{\mu(q)}{\phi(q) z}
\Bigr\vert^{2}
\ \dx \eta
\ll 
N \log N
\]
and
\[
 \sum_{q = 1}^Q 
    \sumast_{a = 1}^q  
    \int_{\xi_{{q,a}}}
\Bigl\vert
\widetilde{S}  \Bigl( \frac aq + \eta \Bigr)  - \frac{\mu(q)}{\phi(q) z}
\Bigr\vert^{2}
\ \dx \eta
\ll 
N \log N.
\]
\end{Lemma} 
\begin{Proof}
By Parseval's theorem and the Prime Number Theorem we have
\[
\int_{-\frac{1}{2}}^{\frac{1}{2}} 
\vert
\widetilde{S}(\alpha) 
\vert^{2}
\ \dx \alpha 
=
 \sum_{m = 1}^{\infty} \Lambda^{2}(m) e^{-2m / N}  
= 
\frac{N}{2} \log N
+
\Odi{N}.
\]
Recalling that the equation at the beginning of page 318
of \cite{LanguascoP1994} implies
\[
\int_{-1/qQ}^{1/qQ}  
\frac{\dx \eta}{\vert z\vert^{2}}  
=
\frac{N}{\pi}\arctan \Bigl(\frac{2\pi N}{qQ} \Bigr)
\]
and using Lemma~2 of Goldston \cite{Goldston1992}, we have
\[
 \sum_{q = 1}^Q 
    \sumast_{a = 1}^q  
    \int_{\xi_{{q,a}}}
\Bigl\vert  \frac{\mu(q)}{\phi(q) z}
\Bigr\vert^{2}
\ \dx \eta
\ll
 \sum_{q = 1}^Q  
  \frac{\mu^2(q)}{\phi(q)}
\int_{-1/qQ}^{1/qQ} 
\frac{\dx \eta}{\vert z\vert^{2}}  
\ll
N \log Q.
\]
The Lemma immediately follows using $Q\leq N$, the relation
$\vert a - b \vert^{2} = \vert a\vert^{2} + \vert b \vert^{2} - 2\Re(a\overline{b})$
and the Cauchy-Schwarz inequality.
\end{Proof}
Let $x\geq 2$ be a real number, $j\geq 0$, $q\geq 1$ be integers and $\chi$ be a Dirichlet character defined $\bmod{q}$.  We define
\begin{equation}
\label{psi-j-def}
\psi_{j}(x ,\chi)
:=
\frac{1}{j!}
\sum_{m=1}^{x} 
(x-m)^{j}
\Lambda(m) \chi(m).
\end{equation}
\begin{Lemma}
\label{explicit-form}
Let $x\geq 2$ be a real number, $j\geq 0$, $q\geq 1$ be integers and $\chi$ be a Dirichlet character defined ${}\bmod{q}$. Assuming that GRH holds
for $L(s, \chi)$ then  
\begin{align*}
\psi_{j}(x ,\chi)
&=
\delta(\chi)
\frac{x^{j+1}}{(j+1)!}
-
\sum_{\rho}
\frac{x^{\rho+ j}}{\rho(\rho + 1) \cdots (\rho + j)}
+
\Odip{j}{x^{j}E(q, x)}+
\Oprimedi{\log x},
\end{align*}
where $\delta(\chi)=1$ if $\chi=\chi_{0} \bmod{q}$ and $0$ otherwise,
$\rho$ runs over the non-trivial zeros of $L(s,\chi)$
and  $E(q, x)$ is defined in \eqref{expl-form-err-term-bis}.
The prime in the last error term means that it is present if and only if  $j=0$.
For $j=0$ the summation over the zeros at the right hand side
should be understood in the symmetric sense.
\end{Lemma}
\begin{Proof}
For $j=0$ this is a classical result, see \emph{e.g.} 
Davenport \cite{Davenport2000}, \S17 and 19.
For $j\geq 1$ it follows by  a standard
Mellin inversion argument using the Ces\`aro kernel defined at
page~142 of Montgomery-Vaughan \cite{MontgomeryV2007}.
We sketch here the proof.

If $q=1$, by eq.~(5.19) of Montgomery-Vaughan \cite{MontgomeryV2007},
we have that
\[
\psi_{j}(x)
=
\psi_{j}(x, \chi_0)
=
-
\frac{1}{2\pi i}
\int_{(c)} 
\frac{\zeta'}{\zeta}(w) 
\frac{x^{w+j}}{w(w + 1) \cdots (w + j)}
 \ \dx w,
\]
where $c>1$ is fixed.
Moving the  line of integration to $\Re(w)=-3/4$, we see that 
the relevant poles are located at the zeros of $\zeta(w)$
and at $w=0,1$. They are all simple poles.
By the Riemann-von Mangoldt formula,
we can choose a large $T$ such that $\vert \Im(\rho)-T \vert\gg (\log T)^{-1}$
for every non-trivial zero $\rho$ of $\zeta(w)$ and hence
there is no harm in moving the integration line to $-3/4$ since
$\vert (\zeta'/\zeta)(w)\vert \ll \log^{2} T$
for every $w=\sigma\pm iT$ with $\sigma \in [-1,2]$, and
$\vert w(w + 1) \cdots (w + j) \vert \gg T^{j+1}$.

By the residue theorem we immediately get
\begin{align}
\notag
\psi_{j}(x)
&=
\frac{x^{j+1}}{(j+1)!}
-
\sum_{\rho}
\frac{x^{\rho+ j}}{\rho(\rho + 1) \cdots (\rho + j)}
-
\frac{\zeta'}{\zeta}(0)
x^{j} \\
\label{zeta-case}
&\qquad
-\frac{1}{2\pi i}
\int_{(-3/4)} 
\frac{\zeta'}{\zeta}(w) 
\frac{x^{w+j}}{w(w + 1) \cdots (w + j)}
\, \dx w.
\end{align}

The vertical integral can be estimated using \eqref{L'/L-bound}
in this special case ($q=1$).
Its contribution is
\begin{equation}
\label{vertical}
\ll
x^{j-3/4}
\int_{-\infty}^{\infty} 
\frac{\log(\vert t \vert+2)}{(1+\vert t \vert)^{j+1}}
\dx t
\ll_{j}
x^{j-3/4}.
\end{equation}
Combining \eqref{zeta-case}-\eqref{vertical} we get the final result
in this case ($j\geq 1$, $q=1$).

Let now $q\geq 2$.
If $\chi$ is the principal character ${}\bmod{q}$ then
\begin{equation}
\label{difference}
\vert
\psi_{j}(x)
-
\psi_{j}(x ,\chi_0)
\vert
\ll
\frac{1}{j!}
\sum_{\substack{{m\leq x}\\ (m,q)>1}}
(x-m)^{j}
\Lambda(m)
\ll
\frac{1}{j!}
x^{j}\log x \log q
\end{equation}
and the result follows using \eqref{zeta-case}-\eqref{vertical}.

Now assume that $\chi\bmod{q}$ is not the principal character ${}\bmod{q}$.
If $\chi\bmod{q}$ were induced by  $\chi_1\bmod{q_1}$, $q_{1} \mid q$,
then, arguing as in \eqref{difference}, we would have
\begin{equation}
\label{chi-difference}
\vert
\psi_{j}(x ,\chi)
-
\psi_{j}(x ,\chi_1)
\vert
\ll
\frac{1}{j!}
x^{j}\log x \log q.
\end{equation}

Now assume that $\chi$ is a primitive character ${}\bmod{q}$.
By eq.~(5.19) of Montgomery-Vaughan \cite{MontgomeryV2007},
we have that
\[
\psi_{j}(x ,\chi)
=
-
\frac{1}{2\pi i}
\int_{(c)} 
\frac{L'}{L}(w,\chi) 
\frac{x^{w+j}}{w(w + 1) \cdots (w + j)}
 \ \dx w,
\]
where $c>1$ is fixed.

We move the  line of integration to $\Re(w)=-3/4$.
The relevant poles are located at the zeros of $L(w,\chi)$
and at $w=0$. They are all simple poles with the unique 
exception of $w=0$ which is a double pole for
$(L'/L)(w,\chi)$ when $\chi$ is even.
By the Riemann-von Mangoldt formula, 
we can choose a large $T$ such that $\vert \Im(\rho)-T \vert\gg (\log (qT))^{-1}$
for every non-trivial zero $\rho$ of $L(s,\chi)$ and hence
there is no harm in moving the integration line to $-3/4$ since 
$\vert (L'/L)(w,\chi)\vert \ll \log^{2} (qT)$
for every $w=\sigma\pm iT$ with $\sigma \in [-1,2]$, and
$\vert w(w + 1) \cdots (w + j) \vert \gg T^{j+1}$. 

By the residue theorem we immediately get
\begin{align}
\notag
\psi_{j}(x ,\chi)
&=
-
\sum_{\rho}
\frac{x^{\rho+ j}}{\rho(\rho + 1) \cdots (\rho + j)}
+ 
C(\chi) \frac{x^{j}}{j!}
\\
\label{general-chi-case}
&\qquad
-\frac{1}{2\pi i}
\int_{(-3/4)} 
\frac{L'}{L}(w,\chi) 
\frac{x^{w+j}}{w(w + 1) \cdots (w + j)}
 \, \dx w,
\end{align}
where $C(\chi)$ is a term that depends only on the character $\chi$.
The vertical integral can be estimated using \eqref{L'/L-bound}
and its contribution is
\begin{equation}
\label{vertical-L}
\ll
x^{j-3/4}
\int_{-\infty}^{\infty} 
\frac{\log(q(1+\vert t \vert))}{(\vert t \vert+2)^{j+1}}
\dx t
\ll_{j}
x^{j-3/4} \log q.
\end{equation}

Finally, we have to deal with the term $C(\chi)$ in
\eqref{general-chi-case}.
If $\chi$ is odd $C(\chi)$ is simply $-(L' / L)(0, \chi) = -b(\chi)$
since $L(0, \chi) \ne 0$.
If $\chi$ is even then $L(w, \chi)$ has a simple zero at $0$ and
therefore $-(L' / L)(w, \chi)x^{w+j}(w(w + 1) \cdots (w + j))^{-1}$ has a double pole at
$w = 0$ with residue $x^{j}C(\chi)/(j!) $ and $C(\chi)=-\log x - b(\chi)$.
The remaining part of the argument runs as at the bottom of Lemma \ref{Linnik-lemma}.
This, together with  \eqref{general-chi-case}-\eqref{vertical-L}, gives the final result
for $j\geq 1$ and a primitive character $\chi\bmod{q}$. 
The general result for $q\geq 2$, $j\geq 1$, follows using \eqref{chi-difference}.
 \end{Proof}

\begin{Lemma}  
\label{I2-lemma}
Let $k\geq 2$ be an integer, 
$N$ be a large integer and $z$ be as in \eqref{z-def}. 
We have 
\begin{align}
\notag
  \int_{-1/2}^{1/2}
&
   W(\chi,\eta,V)
    V(\eta)^{k-1}e(-n \eta) \, \dx \eta   
=
e^{-n / N} \Bigl( \psi_{k-2}(n ,\chi) - \delta(\chi) \frac{n^{k-1}}{(k-1)!}\Bigr)
  +
\Odip{k}{n^{k-2}},
\end{align}
where $\psi_{k-2}(n,\chi)$ is defined in \eqref{psi-j-def}.
\end{Lemma} 
\begin{Proof}
We have
\begin{align*}
  \int_{-1/2}^{1/2}
&
   W(\chi,\eta,V)
    V(\eta)^{k-1}e(-n \eta) \, \dx \eta   \\
&=
 \sum_{m_1 = 1}^{\infty} 
   \sum_{m_2 = 1}^{\infty}
\dotsm
   \sum_{m_k = 1}^{\infty} 
\left( \Lambda(m_1)\chi(m_1) - \delta(\chi) \right) 
e^{- (\sum_{i=1}^k m_i) / N}
   \int_{-1/2}^{1/2} 
   e\Bigl(\bigl(\sum_{i=1}^k m_i - n\bigr) \eta\Bigr) \, \dx \eta \\
\notag
 &=
 \sum_{m_1 = 1}^{\infty} 
   \sum_{m_2 = 1}^{\infty}
\dotsm
   \sum_{m_k = 1}^{\infty} 
\left( \Lambda(m_1)\chi(m_1) - \delta(\chi) \right) e^{- (\sum_{i=1}^k m_i) / N}
  \begin{cases}
     1 & \text{if $\sum_{i=1}^k m_i = n$} \\
     0 & \text{otherwise}
   \end{cases} \\
\label{cancellation}
  &=
    e^{-n / N}
    \sum_{m_1 = 1}^{n - 1} \left( \Lambda(m_1)\chi(m_1) - \delta(\chi) \right)
\binom{n-1-m_1}{k-2}\\
  &=
   \frac{e^{-n / N}}{(k-2)!}
    \sum_{m_1 = 1}^{n - 1} (n-1-m_1)^{k-2} \left( \Lambda(m_1)\chi(m_1) - \delta(\chi) \right) \\
  &\qquad\qquad+
\Odip{k}{n^{k-3} \sum_{m_1 = 1}^{n - 1} (\Lambda(m_1)+1)
}
\\
  &=
    e^{-n / N} \Bigl(\psi_{k-2}(n - 1,\chi) - \delta(\chi) \frac{(n - 1)^{k-1}}{(k-1)!} \Bigr)
+
\Odip{k}{n^{k-2}},
\end{align*}
since the condition $\sum_{i=1}^k m_i= n$ implies that the variables are all $< n$.
Now
\[
\psi_{k-2}(n,\chi) = \psi_{k-2}(n - 1,\chi) +
\Odip{k}{n^{k-3} \sum_{m = 1}^{n - 1} \Lambda(m)}
  =
  \psi_{k-2}(n - 1,\chi)
  +
  \Odip{k}{n^{k-2}}
\]
so that
\[
  e^{-n / N} \Bigl( \psi_{k-2}(n - 1,\chi) - \delta(\chi) \frac{(n - 1)^{k-1}}{(k-1)!} \Bigr)
  =
   e^{-n / N} \Bigl( \psi_{k-2}(n ,\chi) - \delta(\chi) \frac{n^{k-1}}{(k-1)!}\Bigr)
  +
\Odip{k}{n^{k-2}}
\]
and Lemma \ref{I2-lemma} follows.
\end{Proof}

The next lemma is a modern version of Lemma 9 of Hardy-Littlewood \cite{HardyL1923}
and should be compared with equation (1.15) of \cite{FriedlanderG1997}.
\begin{Lemma}
\label{max-estimate-lemma}
Assume GRH,   $1\leq q \leq Q$, $Q\leq N$, $\eta\in \xi_{q,a}$ and let $z$ be as in \eqref{z-def}. Then
\[
\left\vert 
\widetilde{S}\left(\frac{a}{q}+\eta\right) -
\frac{\mu(q)}{\phi(q) z}
\right\vert 
\ll 
\left(
N(q\vert \eta\vert )^{1/2} +(qN)^{1/2}
\right) 
\log N.
\]
\end{Lemma}

\begin{Proof}
By \eqref{approx-with-z-2}, Lemma \ref{Linnik-lemma}
and straightforward computations, we have
\begin{equation}
\label{L1-estim1}
\left\vert 
\widetilde{S}\left(\frac{a}{q}+\eta\right) - \frac{\mu(q)}{\phi(q) z}
\right\vert 
\ll
\frac{q^{1/2}}{\phi(q)}
\Bigl(
\sum_{\chi\bmod{q}}
\vert  
\sum_{\rho}z^{-\rho}\Gamma(\rho)
\vert 
\Bigr) 
+
q^{1/2}\log^2 (qN).
\end{equation}
Since $z^{-\rho} =\vert z\vert ^{-\rho} \exp\left(-i\rho\arctan2\pi N\eta\right)$,
by Stirling's formula we have 
\[
\sum_{\rho}z^{-\rho}\Gamma(\rho) 
\ll
\sum_{\rho}\vert z\vert ^{-1/2}
\exp
\left(
\gamma\arctan2\pi N\eta 
-
\frac{\pi}{2}\vert \gamma\vert 
\right).
\]
If $\gamma\eta \leq 0$ or $\vert \eta\vert \leq 1/N$ we obtain
\begin{equation}
\label{L1-estim2}
\sum_{\rho}z^{-\rho}\Gamma(\rho) 
\ll 
N^{1/2}\log (q+1) ,
\end{equation} 
where, in the first case, $\rho$ runs over the zeros with 
$\gamma\eta \leq
0$.

We can consider only the case $\gamma\eta > 0$ and $\vert \eta\vert > 1/N$.
So we get
\begin{align*} 
\sum_{\rho}z^{-\rho}\Gamma(\rho) 
&\ll
\sum_{\gamma>0}\vert z\vert ^{-1/2}
\exp\Bigl(-\gamma\arctan(\frac{1}{2\pi N\eta})\Bigr) 
\\
&+ \sum_{\gamma<0}\vert z\vert ^{-1/2}
\exp\Bigl(
-\vert \gamma\vert \arctan(\frac{1}{2\pi N\eta})
\Bigr) .
\end{align*} 
We investigate only the case $\gamma>0$ since the other one is
similar. We split $\sum_{\gamma>0}$ according to the cases $\gamma>1$ and
$\gamma\leq 1$ and we denote the first sum as $\sum_1$ and the second
one as $\sum_2$.
Hence, using \eqref{z-estim}, we have
\begin{align}
\notag
\sideset{}{_1}\sum &\ll \vert z\vert ^{-1/2} 
\sum_{m=1}^\infty \log(q(m+1)) \exp
\Bigl(-m\arctan(\frac{1}{2\pi N\eta})\Bigr) 
\\
&
\label{S1-estim} 
 \ll 
 \vert 
 z\vert ^{-1/2} N\eta\log(qN)
 \ll
N\vert \eta\vert ^{1/2}\log(qN).
\end{align} 
Arguing analogously we obtain
\begin{align}
\notag 
\sideset{}{_2}\sum &\ll \vert z\vert ^{-1/2} 
\sum_{0<\gamma\leq 1} 
\exp
\Bigl(-\gamma\arctan(\frac{1}{2\pi N\eta})\Bigr) 
\\
\label{S2-estim} 
&
\ll \vert z\vert ^{-1/2} \log(q+1)
\ll \vert \eta\vert ^{-1/2}\log(qN).
\end{align} 
Lemma \ref{max-estimate-lemma} now follows inserting 
\eqref{L1-estim2}-\eqref{S2-estim} in \eqref{L1-estim1}.
\end{Proof}

Our next lemma concerns the mean-square of the quantity studied in Lemma 1
and it should be considered as a sharper version of equation (7.15) of
Friedlander-Gold\-ston \cite{FriedlanderG1997}. Its proof follows the
argument in Theorem 1 of Languasco-Perelli \cite{LanguascoP1994}:
see also section 5 of \cite{LanguascoP1996}. We insert here just the
relevant changes.
\begin{Lemma}
\label{LP-Lemma}
Assume GRH, let $z$ be as in \eqref{z-def},
$1\leq q \leq Q$ and $Q < N^{1/2}$. Then
\[
 \sumast_{a = 1}^q
\int_{-1/qQ}^{1/qQ}
\Bigl\vert 
\widetilde{S}\Bigl( \frac{a}{q}+\eta \Bigr) - \frac{\mu(q)}{\phi(q) z}
\Bigr\vert ^2 
\dx \eta
\ll
\frac{N}{Q} \log^2 N .
\]
\end{Lemma}
\begin{Proof}
Assuming GRH, by \eqref{approx-with-z-2}-\eqref{expl-form-err-term-bis}   we have
\[
\widetilde{S}\left(\frac{a}{q}+\eta\right) - \frac{\mu(q)}{\phi(q) z} 
=
-
 \frac{1}{\phi(q)}
  \sum_{\chi \bmod q}
    \chi(a) \tau(\overline{\chi})
  \sum_{\rho}z^{-\rho}\Gamma(\rho) 
  +
  \Odi{q^{1/2}\log^2(q N)}
\]
where $\rho=1/2+i\gamma$ runs over the non-trivial zeros of
$L(s,\chi)$ and $\chi\bmod{q}$ is a Dirichlet character.
By the character orthogonality and the previous equation we have
\begin{align}
\notag
 \sumast_{a = 1}^q
\int_{-1/qQ}^{1/qQ}
&
\Bigl\vert 
\widetilde{S}\Bigl( \frac{a}{q}+\eta \Bigr) - \frac{\mu(q)}{\phi(q) z}
\Bigr\vert ^2 
\dx \eta 
\\
\label{ortogonality}
&
\ll 
\frac{q}{\phi(q)} \sum_{\chi\bmod{q}}
\int_{-1/qQ}^{1/qQ}
\Bigl\vert 
\sum_{\rho}z^{-\rho}\Gamma(\rho) 
\Bigr\vert ^2 
\dx \eta 
+ 
\frac{q\log^{4}(q N)}{Q}.
\end{align}
Since $z^{-\rho} =\vert z \vert ^{-\rho} \exp\left(-i\rho\arctan2\pi N\eta\right)$,
by the Stirling formula we have that
\[
\sum_{\rho}
z^{-\rho}\Gamma(\rho)
\ll
\sum_{\rho}
\vert z \vert ^{-1/2}
\exp
\left(
\gamma\arctan2\pi N\eta - \frac{\pi}{2} \vert \gamma\vert 
\right) .
\]
If $\gamma\eta \leq 0$ or $\vert \eta\vert \leq 1/N$ we get,
by \eqref{z-estim}, that
\[
\sum_{\rho}
z^{-\rho}\Gamma(\rho) 
\ll 
N^{1/2} \log (q+1) ,
\]
where $c_{1}>0$ is an absolute constant and, in the first case, 
$\rho$ runs over the zeros with $\gamma\eta \leq 0$.

Let  $\xi=1/(qQ)$. From $1\leq q \leq Q$ and $Q < N^{1/2}$, we have 
$\xi > 1/N$ and  we obtain
\begin{align}
\notag
\int_{-\xi}^{\xi}
\Bigl\vert
\sum_{\rho}
z^{-\rho}\Gamma(\rho)
\Bigr\vert^{2}
\dx\eta 
&
\ll 
\int_{1/N}^{\xi}
\Big\vert 
\sum_{\gamma>0}
z^{-\rho}\Gamma(\rho)\Big\vert ^2 
\dx\eta 
+
\int_{-\xi}^{-1/N}
\Bigl\vert 
\sum_{\gamma<0}
z^{-\rho}\Gamma(\rho)
\Bigr\vert ^2 
\dx\eta \\
\label{coda-alti-bassi}
&
+ 
N \xi \log^{2} (q+1)  .
\end{align}
We will treat only the first integral on the right hand
side of \eqref{coda-alti-bassi}, the second being completely similar. 
Clearly
\begin{equation}
\label{split-int}
\int_{1/N}^{\xi}
\Big\vert 
\sum_{\gamma>0}
z^{-\rho}\Gamma(\rho)
\Big\vert ^2 \dx\eta 
= 
\sum_{k=1}^K
\int_\tau^{2\tau} 
\Big\vert 
\sum_{\gamma>0}
z^{-\rho}\Gamma(\rho)\Big\vert ^2 \dx\eta 
+
\Odi{1}
\end{equation}
where $\tau=\tau_k=\xi/2^k$, $1/N \leq \tau \leq \xi/2$  
and $K$ is a suitable integer satisfying $K=\Odi{\log N}$. 
We can now proceed exactly as at page 312-314 of \cite{LanguascoP1994}.
We obtain
\begin{equation}
\label{int-estim}
\int_{\tau}^{2\tau} 
\Big\vert 
\sum_{\gamma>0}
z^{-\rho}\Gamma(\rho)
\Big\vert ^2
\dx\eta 
\ll 
\sum_{\gamma_1>0}
\sum_{\gamma_2>0}
\frac{1+\Bigl(\frac{\gamma_1+\gamma_2}{N\tau}\Bigr)^2}
{1+\vert \gamma_1-\gamma_2\vert ^2}
\exp\Bigl(
-c\Bigl(\frac{\gamma_1+\gamma_2}{N\tau}
\Bigr)\Bigr) .
\end{equation}
But
\[
\Bigl\{
1+
\Bigl(
\frac{\gamma_1+\gamma_2}{N\tau}
\Bigr)^2
\Bigr\}
\exp\Bigl(
-c_{2}\Bigl(\frac{\gamma_1+\gamma_2}{N\tau}
\Bigr)
\Bigr)
 \ll
\exp\Bigl(
-\frac{c_{2}}{2}\frac{\gamma_1}{N\tau}
\Bigr),
\]
hence the right hand side of \eqref{int-estim} becomes
\begin{equation}
\label{int-estim1}
\ll
\sum_{\gamma_1>0}
\exp\Bigl(
-\frac{c_{2}}{2}\frac{\gamma_1}{N\tau}
\Bigr)
\sum_{\gamma_2>0}
\frac{1}{1+\vert \gamma_1-\gamma_2\vert ^2}.
\end{equation}
Since the number of zeros $\rho_2=1/2+i\gamma_2$ with 
$m\leq \vert \gamma_1-\gamma_2\vert \leq m+1$
is $\Odi{\log (q(m+\vert \gamma_1\vert ))}$,
we immediately get 
\begin{align*}
\notag
\sum_{\gamma_1>0}
\exp\Bigl(
-\frac{c_{2}}{2}\frac{\gamma_1}{N\tau}
\Bigr)
\sum_{\gamma_2>0}
\frac{1}{1+\vert \gamma_1-\gamma_2\vert ^2}
\ll
\int _{0}^{\infty} 
(\log^{2} (qt) )
\exp\Bigl(
-\frac{c_{2}}{2}\frac{t}{N\tau}
\Bigr) 
\dx t.
\end{align*}

The  function $(\log^{2} (qt) )
\exp(
-\frac{c_{2}}{4}\frac{t}{N\tau}
) $
has a maximum 
attained at $t_{0}$ such that $t_{0}\log (qt_{0}) $ $= 8N\tau/c_{2}$.
Hence the right hand side of the previous equation is
\begin{align}
\notag
&
\ll
\int _{0}^{1/q} 
(\log^{2} (qt) )
\exp\Bigl(
-\frac{c_{2}}{2}\frac{t}{N\tau}
\Bigr) 
\dx t 
+
\int _{1/q}^{t_{0}} 
(\log^{2} (qt) )
\exp\Bigl(
-\frac{c_{2}}{2}\frac{t}{N\tau}
\Bigr) 
\dx t
\\
\notag
&\qquad
+
\int _{t_{0}}^{\infty} 
(\log^{2} (qt) )
\exp\Bigl(
-\frac{c_{2}}{2}\frac{t}{N\tau}
\Bigr) 
\dx t
\\
\notag
&
\ll
t_{0}(\log^{2} (qt_{0}) )
+
N\tau(\log^{2} (qt_{0})  )
\exp\Bigl(
-\frac{c_{2}}{2}\frac{t_{0}}{N\tau}
\Bigr) 
\\
\label{int-estim1-alto}
&
\ll
N\tau(\log^{2} (qt_{0})  )
\exp\Bigl(
-\frac{c_{2}}{2}\frac{t_{0}}{N\tau}
\Bigr) 
\ll
N\tau \log^{2} (qN).
\end{align}
Hence, inserting 
\eqref{int-estim1-alto} into 
\eqref{int-estim}-\eqref{int-estim1}, we get
\begin{equation}
\label{int-estim2-alto}
\int_{\tau}^{2\tau} 
\Big\vert 
\sum_{\gamma>0}
z^{-\rho}\Gamma(\rho)
\Big\vert ^2
\dx\eta 
\ll
N\tau \log^{2} (qN).
\end{equation}
Inserting now
\eqref{int-estim2-alto} 
into \eqref{coda-alti-bassi}-\eqref{split-int}
we get
\[
\int_{-1/qQ}^{1/qQ}
\Bigl\vert 
\sum_{\rho}z^{-\rho}\Gamma(\rho) 
\Bigr\vert ^2 
\dx \eta 
\ll 
\frac{N}{qQ} \log^{2} (qN)
\]
and hence, by \eqref{ortogonality},  Lemma \ref{LP-Lemma}  follows. 
\end{Proof}

The next lemma will be useful in the computation of the main term in
Theorem~\ref{Explicit-formula-Theorem}.
We insert the proof already contained in Languasco-Perelli \cite{LanguascoP1994}
for $k=2$ and in Languasco \cite{Languasco2000a} for $k\geq 3$.
\begin{Lemma}
\label{main-term-lemma}
Let $k\geq 2$ be an integer, $z$ be as in \eqref{z-def} and $Q\leq N^{1/2}/2$. 
Then, uniformly for $1\leq \ell\leq N$, we have
\[
\int_{\xi_{q,a}} \frac{e(-\ell\eta)}{z^k} \, \dx \eta 
=
e^{-\ell/N} \frac{\ell^{k-1}}{(k-1)!} 
+
\Odi{(qQ)^{k-1}}.
\]
\end{Lemma}
\begin{Proof}
Let $T\geq 1/2$. Using \eqref{z-def} we get 
\begin{align}
\notag 
\int_{\xi_{q,a}} \frac{e(-\ell\eta)}{z^k} \, \dx \eta 
&=
\int_{-T}^{T} \frac{e(-\ell\eta)}{z^k} \, \dx \eta  
+ 
\Odi{
\int_{\frac{1}{2qQ}}^{T} \frac{\dx \eta}{\vert z\vert ^k}+
\int_{-T}^{-\frac{1}{2qQ}} \frac{\dx \eta}{\vert z\vert ^k}}
\\
\label{L5-estim1}
&
=
\int_{-T}^{T} \frac{e(-\ell\eta)}{z^k} \, \dx \eta + 
\Odi{(qQ)^{k-1}}.
\end{align} 
Using the variable $z = N^{-1} -2\pi i \eta$ in place of $\eta$, we have
\begin{equation}
\label{substit}
\int_{-T}^T \frac{e(-\ell\eta)}{z^k} \, \dx \eta 
= 
\frac{e^{-\ell/N}}{2\pi i}
\int_{\frac{1}{N}-2\pi iT}^{\frac{1}{N}+2\pi iT} 
\frac{\exp(\ell z)}{z^k} \, \dx z.
\end{equation}

Let $\Gamma$ denote the left half of the circle
$\left\vert z- N^{-1}\right\vert =2\pi T$.
By the residue theorem we obtain 
\begin{align}
\notag
\frac{e^{-\ell /N}}{2\pi i} 
\int_{\frac{1}{N}-2\pi iT}^{\frac{1}{N}+2\pi iT} 
\frac{\exp(\ell z)}{z^k} \, \dx z 
&
= 
e^{-\ell /N}\frac{\ell ^{k-1}}{(k-1)!} 
+ 
\frac{e^{-\ell /N}}{2\pi i}\int_{\Gamma} \frac{\exp(\ell z)}{z^k} \, \dx z 
\\
\label{residue}
&= 
e^{-\ell /N}\frac{\ell ^{k-1}}{(k-1)!} 
+ 
\Odi{\frac{1}{T^{k-1}}} .
\end{align}
Lemma \ref{main-term-lemma} now follows from
\eqref{L5-estim1}-\eqref{residue} letting $T\to \infty$.
\end{Proof}

Lemma \ref{aux-lemma} below follows inserting 
Lemma \ref{LP-Lemma} and \eqref{z-def} in the
body of the proof of Lemma 5 of Friedlander-Goldston \cite{FriedlanderG1997}.

\begin{Lemma}
\label{aux-lemma}
Assume GRH and let $z$ be as in \eqref{z-def}. 
Then, for any real $c>0$, we have
\[
\sum_{q=1}^Q
\sumast_{a=1}^q
\frac{1}{\phi(q)^c}
\int_{-1/qQ}^{1/qQ} 
\Bigl \vert
\widetilde{S}\Bigl( \frac{a}{q}+\eta \Bigr)
-
\frac{\mu(q)}{\phi(q)}\frac{1}{z}
\Bigr \vert^2
\Bigl \vert
\frac{\mu(q)}{\phi(q)}\frac{1}{z}
\Bigr \vert^2 
\, \dx \eta
\ll 
N^2 \log^2 N 
\]
and, for $c=0$, the same result holds replacing  $ \log^2 N$
with $ \log^3 N $.
\end{Lemma}

Let now 
\begin{equation}
\label{S-star-def}
  S^*(Q)
  =
  \max_{q \le Q}
    \max_{(a, q) = 1}
      \max_{\eta \in \xi_{q,a}}
        \bigl\vert  \widetilde{R}(\eta; q, a, z) \bigr\vert .
\end{equation}

We have
\begin{Lemma}
\label{mixed-estims}
Let $m\geq 2$ and $z$ be as in \eqref{z-def}. Then
\begin{equation}
\notag
\sum_{q=1}^Q\sumast_{a = 1}^q
\int_{-1/qQ}^{1/qQ} 
\left \vert 
\widetilde{S}\Bigl( \frac{a}{q}+\eta \Bigr) -\frac{\mu(q)}{\phi(q) z}
\right \vert  ^m 
\dx \eta
\ll 
(S^*(Q))^{m-2}N\log N.
\end{equation} 
Assuming GRH we have
\begin{align*}
\sum_{q=1}^Q\sumast_{a = 1}^q
\int_{-1/qQ}^{1/qQ} 
\left \vert 
\widetilde{S}\Bigl( \frac{a}{q}+\eta \Bigr) -\frac{\mu(q)}{\phi(q) z}
\right \vert ^m 
\left \vert 
\frac{\mu(q)}{\phi(q) z}
\right \vert 
\dx \eta
& \ll 
(S^*(Q))^{m-2}N^{3/2}\log^2 N,
\\
\sum_{q=1}^Q\sumast_{a = 1}^q
\int_{-1/qQ}^{1/qQ} 
\left \vert 
\widetilde{S}\Bigl( \frac{a}{q}+\eta \Bigr) -\frac{\mu(q)}{\phi(q) z}
\right \vert ^m 
\left \vert 
\frac{\mu(q)}{\phi(q) z}
\right \vert ^2 
\dx \eta
&
\ll 
(S^*(Q))^{m-2}N^{2}\log^3 N,
\\
\sum_{q=1}^Q\sumast_{a = 1}^q
\int_{-1/qQ}^{1/qQ} 
\left \vert 
\widetilde{S}\Bigl( \frac{a}{q}+\eta \Bigr) -\frac{\mu(q)}{\phi(q) z}
\right \vert 
\left \vert 
\frac{\mu(q)}{\phi(q) z}
\right \vert ^2 
\dx \eta
&
\ll 
N^{3/2}\log^2 N,
\end{align*} 
and, for $r\geq 3$,
\begin{align*}
\sum_{q=1}^Q\sumast_{a = 1}^q
\int_{-1/qQ}^{1/qQ} 
\left \vert 
\widetilde{S}\Bigl( \frac{a}{q}+\eta \Bigr) -\frac{\mu(q)}{\phi(q) z}
\right \vert  ^m 
\left \vert 
\frac{\mu(q)}{\phi(q) z}
\right \vert ^r 
\dx \eta
&
\ll 
(S^*(Q))^{m-2}N^{r}\log^2 N,
\\
\sum_{q=1}^Q\sumast_{a = 1}^q
\int_{-1/qQ}^{1/qQ} 
\left \vert 
\widetilde{S}\Bigl( \frac{a}{q}+\eta \Bigr) -\frac{\mu(q)}{\phi(q) z}\
\right \vert 
\left \vert 
 \frac{\mu(q)}{\phi(q) z}
 \right \vert^r 
 \dx \eta
&
\ll 
N^{r-1/2}(\log N)^{3/2}.
\end{align*} 
\end{Lemma}
The proof of Lemma \ref{mixed-estims} follows using  
Lemmas \ref{incond-mean-square}, 
\ref{LP-Lemma}, \ref{aux-lemma}
arguing as in Lemma 3 of Friedlander-Goldston \cite{FriedlanderG1997}.

\section{Proof of the main result}

We consider the usual Farey dissection of level $Q$ of the unit
interval as in \eqref{Farey-arc}.
By \eqref{approx-with-V-1}, we have
\begin{align}
\notag
  e^{-n / N}
  &R_k(n)
  =
  \int_0^1
    \widetilde{S}(\alpha)^k e(-n \alpha) \, \dx \alpha \\
    \notag
  &=
  \sum_{q = 1}^Q
    \Bigl( \frac{\mu(q)}{\phi(q)} \Bigr)^k
    \sumast_{a = 1}^q e \Bigl( -n \frac aq \Bigr)
      \int_{\xi_{q,a}} V(\eta)^{k}e(-n \eta) \, \dx \eta \\
      \notag
  &\quad+
  k
  \sum_{q = 1}^Q
    \Bigl( \frac{\mu(q)}{\phi(q)} \Bigr)^{k - 1}
    \sumast_{a = 1}^q e \Bigl( -n \frac aq \Bigr)
      \int_{\xi_{q,a}}
        \widetilde{R}(\eta; q, a, V) V(\eta)^{k-1} e(-n \eta) \, \dx \eta \\
        \notag
  &\quad+
  \sum_{m = 2}^k
    \binom km
  \sum_{q = 1}^Q
    \Bigl( \frac{\mu(q)}{\phi(q)} \Bigr)^{k - m}
    \sumast_{a = 1}^q e \Bigl( -n \frac aq \Bigr)
      \int_{\xi_{q,a}}
        \widetilde{R}(\eta; q, a, V)^m V(\eta)^{k-m}e(-n \eta) \, \dx \eta \\
  \label{terms-def}
  &=
  M_0(k) + k M_1(k)
  +
  \sum_{m = 2}^k
    \binom km M_m(k),
\end{align}
say.

\subsection{Main term $M_0(k)$}
By Lemma \ref{V-behaviour} we can write
\[
  \int_{\xi_{q,a}} V(\eta)^{k} e(-n \eta) \, \dx \eta
=
  \int_{\xi_{q,a}} \frac{e(-n \eta)}{z^k} \, \dx \eta
+
\Odi{
  \int_{-1/(qQ)}^{1/(qQ)} \frac{\dx \eta}{\vert z \vert^{k-1}}
}
\]
and, using \eqref{z-estim}, the error term in the previous equation is
\[
\ll
\int_{-1/N}^{1/N} N^{k-1} \, \dx \eta
+
\int_{1/N}^{1/(qQ)} \frac{\dx \eta}{\vert \eta \vert^{k-1}} 
  \ll_k
  \begin{cases}
    \log(N / q Q) & \textrm{if $k = 2$,} \\
    N^{k - 2}                 & \textrm{if $k > 2$.}
  \end{cases}
\]
Combining the previous two equations with Lemma \ref{main-term-lemma},
for $k\geq 3$ we have
\[
  \int_{\xi_{q,a}} V(\eta)^{k} e(-n \eta) \, \dx \eta
 =
  e^{-\ell / N} \frac{n^{k - 1}}{(k - 1)!}
  +
  \Odip{k}{(q Q)^{k - 1}+N^{k-2}}
\]
uniformly for $1 \le n \le N$ and $Q \le N^{1 / 2}/2$.
Recalling the definition for the Ramanujan sum $c_q$ in \eqref{cq-def}
and for the singular series $\singseries_k(n)$ in \eqref{singseries-def},
we get
\begin{align}
\notag
  M_0(k)
  &=
   e^{-n / N}
  \sum_{q = 1}^Q
    \Bigl( \frac{\mu(q)}{\phi(q)} \Bigr)^k
    c_q(-n) \frac{n^{k - 1}}{(k - 1)!}
  +
  \Odi{ \sum_{q=1}^Q \frac{\vert c_q(-n) \mu(q) \vert }{\phi(q)^k} ((q Q)^{k - 1}+ N^{k-2})} \\
  \notag
  &=
  e^{-n / N} \frac{n^{k - 1}}{(k - 1)!}
  \singseries_k(n)
  +
  \Odi{ n^{k - 1} \sum_{q > Q} \frac{\mu^2(q)}{\phi(q)^k} \vert c_q(-n)\vert  } \\
  \notag
  &\qquad+
  \Odi{Q^{k - 1} \sum_{q = 1}^Q \mu^2(q) \Bigl( \frac q{\phi(q)} \Bigr)^{k - 1} } 
  + 
  \Odi{N^{k - 2} \sum_{q = 1}^Q  \frac 1{\phi(q)^{k - 1}} } \\
  \label{M0-eval}
  &=
  e^{-n / N} \frac{n^{k - 1}}{(k - 1)!}
  \singseries_k(n)
  +
  \Odi{Q^k + n^{k - 1 + \eps} Q^{1 - k} + N^{k - 2}},
\end{align}
since $k \ge 3$ and
\begin{align*}
 \sum_{q > Q} \frac{\mu^2(q)}{\phi(q)^k} \vert c_q(-n)\vert 
 &\leq
  \sum_{q > Q} \frac{\mu^2(q)}{\phi(q)^k} \sum_{\substack{d\mid n\\ d\mid q}} d
   \leq
   \sum_{d\mid n}\frac{d\mu^2(d)}{\phi(d)^k}
  \sum_{q'> Q/d} \frac{\mu^2(q')}{\phi(q')^k} \\
 & \ll
  Q^{1-k}   \sum_{d\mid n}\frac{d^{k}\mu^2(d)}{\phi(d)^k}
    \ll
  Q^{1-k} n^\eps,
\end{align*}
using also Lemma~2 of Goldston \cite{Goldston1992}.
Then for $n=N$, $k \ge 3$ and $Q = N^{1 / 2}/2$ the error terms 
are under control, since we have to compare them with the 
order of magnitude of the secondary main term which is 
$\approx N^{k - 3/2}$.

\subsection{Secondary main term $M_1(k)$}
\label{secondary-main-term}

Equation \eqref{terms-def} implies that
\[
  M_1(k)
  =
  \sum_{q = 1}^Q
    \Bigl( \frac{\mu(q)}{\phi(q)} \Bigr)^{k - 1}
    \sumast_{a = 1}^q e \Bigl( -n \frac aq \Bigr)
      \int_{\xi_{q,a}}
        \widetilde{R}(\eta; q, a, V) V(\eta)^{k-1} e(-n \eta) \, \dx \eta.
\]
We set $\theta_{q,a} = (-1/2, 1/2) \setminus \xi_{q,a}$ so
that $M_1(k) = A - B$, say, where
\begin{align*}
  A
  &:=
  \sum_{q = 1}^Q
    \Bigl( \frac{\mu(q)}{\phi(q)} \Bigr)^{k - 1}
    \sumast_{a = 1}^q e \Bigl( -n \frac aq \Bigr)
      \int_{-1/2}^{1/2}
        \widetilde{R}(\eta; q, a, V) V(\eta)^{k-1} e(-n \eta) \, \dx \eta, \\
  B
  &:=
  \sum_{q = 1}^Q
    \Bigl( \frac{\mu(q)}{\phi(q)} \Bigr)^{k - 1}
    \sumast_{a = 1}^q e \Bigl( -n \frac aq \Bigr)
      \int_{\theta_{q,a}}
        \widetilde{R}(\eta; q, a, V) V(\eta)^{k-1} e(-n \eta) \, \dx \eta.
\end{align*}
In order to estimate $B$, we first remark that
\[ 
    |\widetilde{R}(\eta; q, a, V) |
    \ll 
    \Bigl|\widetilde{S}\Bigl( \frac aq + \eta \Bigr) \Bigr|
  +
  \frac{\mu^2(q)}{\phi(q)} 
    |V(\eta)|
    \ll
      \sum_{n \ge 1}
      \Lambda(n)  e^{-n / N}  
      +
        \frac{N}{\phi(q)}  
        \ll N
\]
by \eqref{V-estim} and the Prime Number Theorem. Hence,
since  $\theta_{q,a} \subset (-1/2, -1/(2qQ)) \cup (1/(2qQ),1/2)$, we obtain
\begin{align*}
  |B|
  & \ll
  N
  \sum_{q=1}^Q
  \frac{\mu^2(q)}{\phi^{k-2}(q)}
  \Bigl(
  \int_{1/(2qQ)}^{1/2}
  +
  \int_{-1/2}^{-1/(2qQ)}
  \Bigr)
  \vert V(\eta) \vert ^{k-1} \dx \eta 
 \ll_k
  N Q^{k-2}
  \sum_{q=1}^Q
  \frac{\mu^2(q)q^{k-2}}{\phi^{k-2}(q)}
  \\
  & \ll_k
  N Q^{k-1},
\end{align*}
by \eqref{V-estim} and  Lemma~2 of Goldston \cite{Goldston1992}. 
We explicitly remark that the usual strategy to estimate $B$ involves
the Cauchy-Schwarz inequality. In this case this would lead to
$\vert B \vert \ll_k (N\log N)^{1/2} Q^{k}$ which is worse than
our estimate for $Q > (N/\log N)^{1/2}$. 
In this case the optimal choice
of $Q$ will be $N^{1/2}/2$, see \S\ref{error-terms} below, and hence
our estimate is slightly sharper.
Summing up,
\begin{align*}
  M_1(k)
  &=
  \sum_{q = 1}^Q
    \Bigl( \frac{\mu(q)}{\phi(q)} \Bigr)^{k - 1}
    \sumast_{a = 1}^q e \Bigl( -n \frac aq \Bigr)
      \int_{-1/2}^{1/2}
        \widetilde{R}(\eta; q, a, V) V(\eta)^{k-1} e(-n \eta) \, \dx \eta 
        \\
  &\qquad
  +
  \Odip{k}{N Q^{k-1}}.
\end{align*}
Inserting the approximation \eqref{approx-with-V-2}, we have
\begin{align*}
  M_1(k)
  &=
  \sum_{q = 1}^Q
    \frac{\mu(q)^{k-1}}{\phi(q)^k}
    \sumast_{a = 1}^q e \Bigl( -n \frac aq \Bigr)
      \sum_{\chi \bmod q} \chi(a) \tau(\overline{\chi})
        \int_{-1/2}^{1/2}
          W(\chi, \eta, V) V(\eta)^{k-1} e(-n \eta) \, \dx \eta \\
  &\qquad+
  \Odip{k}{N Q^{k-1}} 
  +
  \Odip{k}{
    \sum_{q = 1}^Q \frac{\mu^2(q)}{\phi(q)^{k-2}} \log^{2} (qN)
    \int_{-1/2}^{1/2} |V(\eta)|^{k-1} \, \dx \eta
   }.
\end{align*}
The inequality \eqref{V-estim} implies that the last integral above is
$\ll_k N^{k - 2}$, and the last error term is
\[
  \ll_k
  N^{k - 2} \log^{2} (QN)
  \sum_{q = 1}^Q \frac{\mu^2(q)}{\phi(q)^{k-2}}
  \ll_k
  N^{k - 2} \log^{2} (QN) f(Q, k),
\]
where
\[
  f(Q, k)
  =
  \begin{cases}
    \log Q &\text{if $k = 3$,} \\
    1      &\text{if $k \ge 4$.}
  \end{cases}
\]
By  
\eqref{c-chi-def},  we can now write
\begin{align*}
  M_1(k)
  &=
\sum_{q = 1}^Q  
  \frac{\mu(q)^{k-1}}{\phi(q)^{k}} 
    \sum_{\chi \bmod q}
      c_{\chi}(-n) \tau(\overline{\chi})
  \int_{-1/2}^{1/2}
   W(\chi,\eta,V)
    V(\eta)^{k-1}e(-n \eta) \, \dx \eta \\
  &\qquad+
  \Odip{k}{ 
N Q^{k-1} + N^{k - 2} \log^{2} (QN) f(Q, k) 
}.
\end{align*}
Assuming GRH holds and
using Lemma \ref{I2-lemma},  we obtain
\begin{align}
\notag
  M_1(k)
  &=
 e^{-n / N} 
\sum_{q = 1}^Q  
  \frac{\mu(q)^{k-1}}{\phi(q)^{k}} 
    \sum_{\chi \bmod q}
      c_{\chi}(-n) \tau(\overline{\chi})
  \Bigl( \psi_{k-2}(n ,\chi) - \delta(\chi) \frac{n^{k-1}}{(k-1)!}\Bigr)
\\
\notag
&\qquad+
  \Odip{k}{ 
N Q^{k-1}
+ 
N^{k - 2} \log^{2} (QN) f(Q, k)
}\\
\label{almost-done}
&\qquad+
  \Odip{k}{ n^{k-2}
\sum_{q = 1}^Q  
  \frac{\mu(q)^{2}}{\phi(q)^{k}} 
    \sum_{\chi \bmod q}
	\vert c_{\chi}(-n) \tau(\overline{\chi}) \vert
	}.
\end{align}
Using Lemma \ref{explicit-form} with $j=k-2$
and remarking that in this case the error term is $\log^2(Qn)$
times the last error term of \eqref{almost-done},
we obtain
\begin{align}
\notag
   M_1(k)
  &=
  - e^{-n / N} 
  \sum_{q = 1}^Q
    \frac{\mu(q)^{k - 1}}{\phi(q)^k}
    \sum_{\chi \bmod q}
      c_{\chi}(-n) \tau(\overline{\chi})
      \sum_{\rho}
        \frac{n^{\rho+k-2}}{\rho(\rho + 1) \cdots (\rho + k - 2)} \\
        \label{M1-eval}
  &\qquad+
  \Odip{k}{ 
N Q^{k-1}
+ 
N^{k - 2} \log^{2} (QN)f(Q, k)
+
n^{k - 2} \log^{2} (Qn) g(Q, k)
},
\end{align}
where, by Lemma~2 of Goldston \cite{Goldston1992},
the last error term in \eqref{almost-done} is
\begin{equation}
\notag
\ll
n^{k-2}
\sum_{q = 1}^Q  
  \frac{\mu(q)^{2}q^{1/2}}{\phi(q)^{k-2}} 
\ll_{k}
n^{k-2}
  \begin{cases}
    Q ^{1/2} & \textrm{if $k =   3$,} \\
    1 & \textrm{if $k \ge 4$}
  \end{cases}
  = 
n^{k-2}
    g(Q, k).
\end{equation}

We remark that the summation over $q$ can be extended 
to all $q$ by inserting a new error term which is 
$\Odi{n^{k-2+\Theta} Q^{3-k+\epsilon}}$,
where  $\Theta  =  \sup_{\chi \bmod q}\{\beta : L(\beta + i\gamma, \chi) = 0\}$,
see p.~296 of \cite{FriedlanderG1997}. 
For $n=N$ and $Q = N^{1 / 2}/ 2$, the previous unconditional estimate becomes admissible for $k \ge 5$ and hence we can ``just'' assume
that GRH holds for every $q\leq N^{1/2}/2$
(the order of magnitude of the secondary main term is
$\approx N^{k - 3/2}$).
Assuming GRH in its ``full strength'' the tail of the singular series
gives a contribution of $\Odi{n^{k/2+\epsilon}}$ which is admissible for $k\ge 4$.

Hence we can finally say, for $n=N$, $k \ge 5$ and $Q = N^{1 / 2}/ 2$, that
the error terms are under control under the assumption of GRH for every 
$q\leq N^{1 / 2}/ 2$.

 \subsection{The error terms}
\label{error-terms}
Essentially, they are estimated as in Languasco \cite{Languasco2000a}.
Assuming the Generalized Riemann Hypothesis, we have,
by Lemma \ref{max-estimate-lemma} and \eqref{S-star-def}, the estimate
\begin{equation}
\label{minorarcs}
  S^*(Q)
  \ll
  \max_{q \le Q}
    \max_{(a, q) = 1}
      \max_{\eta \in \xi_{q,a}}
        \Bigl[ \log(q N) \bigl( N \sqrt{q \vert \eta\vert } + \sqrt{q N} \bigr)
        \Bigr]
  \ll
  \Bigl(
    \frac{N}{\sqrt{Q}} + \sqrt{Q N}
  \Bigr)
  \log(Q N).
\end{equation}
 
The optimal $Q$ in \eqref{minorarcs} is $Q=N^{1/2}/2$ and in this case we get
\begin{equation}
\label{optimal-Q}
  S^*(Q)
  \ll
 N^{3/4} \log N.
\end{equation}

Using this notation and  Lemma \ref{mixed-estims}, 
we can write the following bounds:

\begin{enumerate}[i)]
\item 
for $m = k$, $m \ge 2$, we unconditionally get
\[
  M_k(k)
  \ll_k
  \bigl( S^*(Q) \bigr)^{k - 2} N \log N;
\]
\item for $m = k - 1$, $m \ge 2$, $k \ge 3$, assuming GRH, 
we obtain
\[
  M_{k-1}(k)
  \ll_k
  \bigl( S^*(Q) \bigr)^{k - 3} N^{3 / 2} \log^2 N;
\]
\item for $m = k - 2$, $m \ge 2$, $k \ge 4$,  assuming GRH, 
we obtain
\[
  M_{k-2}(k)
  \ll_k
  \bigl( S^*(Q) \bigr)^{k - 4} N^2 \log^3 N;
\]
\item for $2 \le m \le k - 3$, $k \ge 5$, assuming GRH, 
we obtain
\[
  M_m(k)
  \ll_k
  \bigl( S^*(Q) \bigr)^{m - 2} N^{k - m} \log^2 N.
\]
\end{enumerate}

\section{Conclusion of the proof}

We restrict our analysis to $k\geq 5$ since the
error terms in \S \ref{secondary-main-term} are under control only in this case.
If $k = 5$, the expected main term has size $N^4$ and error terms $\Odi{N^{3}}$, the secondary main term has expected size $N^{7/2}$ and $M_1(5)$ has an error term $\Odi{N^{3} \log^{2} N}$.
Moreover, using \eqref{optimal-Q}, we get
\begin{align*}
  M_2(5)
  &\ll
  N^3 \log^2 N,
  &
  M_3(5)
  &\ll
  N^{7 / 4} \log^4 N, \\
  M_4(5)
  &\ll
  N^3 \log^4 N,
  &
  M_5(5)
  &\ll
  N^{13 / 4} \log^4 N
\end{align*}
and hence the global error term in this case is $N^{13 / 4} \log^4 N$.

If $k \geq 6$, the expected main term has size $N^{k-1}$ and
error terms $\Odi{N^{k-2}}$,
the secondary main term has expected size $N^{k-3/2}$ and
$M_1(k)$ has an error term $\Odi{N^{k-2} \log^2 N}$.
Moreover, again by \eqref{optimal-Q}, we obtain
\begin{align*}
  M_m(k)
  &\ll_{k}
     N^{ k - 3 / 2 - m / 4}  \log^m N
     &\textrm{for}\quad  & 2 \le m \le k - 3,\\
  M_{k-2}(k)
  &\ll_{k}
    N^{(3 / 4) k - 1}  \log^{k-1} N,
  &M_{k-1}(k)
  &\ll_{k}
  N^{(3 / 4) k - 3/4}  \log^{k-1} N,\\
  M_k(k)
  &\ll_{k}
  N^{(3 / 4) k - 1/2} \log^{k-1} N.
\end{align*}
The maximum for $M_m(k)$ is attained at $m=2$ and is $\ll_{k}N^{k-2} \log^2 N$.
Hence, for $k=6$ the global upper bound is $\ll N^{4} \log^5 N$
while for $k\geq 7$ it is $\ll_{k} N^{k-2} \log^2 N$.

Combining the previous remarks with \eqref{terms-def}-\eqref{M0-eval} and \eqref{M1-eval}, the Theorem follows.

\vskip 1cm
\noindent
Alessandro Languasco, Dipartimento di Matematica Pura e Applicata, Universit\`a
di Padova, Via Trieste 63, 35121 Padova, Italy.\\
\emph{email:} languasco@math.unipd.it

\medskip
\noindent
Alessandro Zaccagnini, Dipartimento di Matematica, Universit\`a di Parma, Parco
Area delle Scienze 53/a, Campus Universitario, 43124 Parma, Italy. \\
\emph{email:} alessandro.zaccagnini@unipr.it

\end{document}